\documentclass[11pt]{amsart}
\usepackage{color}
\usepackage{amsmath}
\usepackage{amssymb}
\usepackage{amsfonts}
\usepackage[latin1]{inputenc}

\DeclareMathSymbol{\subsetneqq}{\mathbin}{AMSb}{36}

\theoremstyle{plain}
 \numberwithin{equation}{section}

\newcommand{\bea}     {\begin{eqnarray}}
\newcommand{\eea}     {\end{eqnarray}}
\newcommand{\beal}     {\begin{align}}
\newcommand{\eeal}     {\end{align}}
\newcommand{\bsp}     {\begin{split}}
\newcommand{\esp}     {\end{split}}
\newcommand{\blll}     {\begin{array}{lll}}
\newcommand{\brcl}     {\begin{array}{rcl}}
\newcommand{\barr}     {\begin{array}}
\newcommand{\earr}     {\end{array}}
\newcommand{\beet}     {\begin{eqnarray*}}
\newcommand{\eec}     {\end{center}}
\newcommand{\bequ}      {\begin{equation}}
\newcommand{\eequ}      {\end{equation}}
\newcommand{\bet}     {\begin{tabular}}
\newcommand{\eet}     {\end{tabular}}
\newcommand{\btab}     {\begin{table}}
\newcommand{\etab}     {\end{table}}
\newcommand{\ben}     {\begin{enumerate}}
\newcommand{\een}     {\end{enumerate}}
\newcommand{\bec}     {\begin{center}}
\newcommand{\beit}     {\begin{itemize}}
\newcommand{\eeit}     {\end{itemize}}

\newtheorem{theorem}{Theorem}[section]

\newtheorem{lemma}{Lemma}[section]
\newtheorem{proposition}{Proposition}[section]
\newtheorem{definition}{Definition}[section]
\newtheorem{remark}{Remark}[section]
\newtheorem{notation}{Notation}

\newcommand{\edoc}{\end{document}}

\newcommand{\rp}{\mathrm{I\kern-.19em R}^{p}}
\newcommand{\rn}{\mathrm{I\kern-.19em R}^{n}}
\newcommand{\rem}{\mathrm{I\kern-.19em R}^{m}}
\newcommand{\rnn}{\mathrm{I\kern-.19em R}^{n\times n}}
\newcommand{\remm}{\mathrm{I\kern-.19em R}^{m\times m}}
\newcommand{\rmn}{\mathrm{I\kern-.19em R}^{m\times n}}
\newcommand{\rnm}{\mathrm{I\kern-.19em R}^{n\times m}}
\newcommand{\rmp}{\mathrm{I\kern-.19em R}^{m\times p}}
\newcommand{\rpn}{\mathrm{I\kern-.19em R}^{p\times n}}
\newcommand{\rpp}{\mathrm{I\kern-.19em R}^{p\times p}}
\newcommand{\rqq}{\mathrm{I\kern-.19em R}^{q\times q}}
\newcommand{\rpq}{\mathrm{I\kern-.19em R}^{p\times q}}
\newcommand{\rppq}{\mathrm{I\kern-.19em R}^{(p+q)\times (p+q)}}
\newcommand{\rnmu}{\mathrm{I\kern-.19em R}^{(n-1)\times (n-1)}}
\newcommand{\rnmd}{\mathrm{I\kern-.19em R}^{(n-2)\times (n-2)}}

\newcommand{\Q}{{\rm \kern.24em \vrule width.02em height1.4ex depth-.05ex \kern-.26em Q}}
\newcommand{\1}{{\rm 1\kern-.25em 1}}
\newcommand{\R}{{\rm I\kern-.25em R}}
\newcommand{\E}{{\rm I\kern-.25em E}}
\newcommand{\N}{{\rm I\kern-.25em N}}
\renewcommand{\H}{{\rm I\kern-.25em H}}

 \input{tcilatex}
 
\begin{document}

\title[Global well-posedness for the MQG equation ]{Global well-posedness for
a Modified 2D dissipative quasi-geostrophic equation with initial
data in the critical Sobolev space $H^1$}
\author{Ramzi May}
\address{Département de Mathématiques, Faculté des Sciences de Bizerte,
Jarzouna 2001 Bizerte Tunisie.} \email{Ramzi.May@fsb.rnu.tn}
\date{3 October 2009}
\maketitle

 \medskip
\begin{abstract}
In this paper, we consider the following modified
quasi-geostrophic equations
\begin{equation}
\partial _{t}\theta +\Lambda ^{\alpha }\theta +u\vec{\nabla}\theta =0,~~~\
u=\Lambda ^{\alpha -1}\mathcal{R}^{\perp }(\theta )  \tag{MQG}
\end{equation}%
wher $\alpha \in ]0,1[$ is a fixed parameter. This equation was
recently introduced by P. Constantin, G. Iyer and J. Wu in
\cite{CIW} as a modification of the classical quasi-geostrophic
equation. In this paper, we prove that for any initial data
$\theta _{\ast }$ in the Sobolev space $H^{1}(\mathbb{R}^{2}),$
the equation (MQG) has a global and smooth solution $\theta $ in $%
C\left( \mathbb{R}^{+},H^{1}\left( \mathbb{R}^{2}\right) \right)
.$
\end{abstract}

\section{Introduction and statement of results}
In this paper, we are concerned with the modified 2D dissipative
quasi-geostrophic equation:%
\begin{equation}
\left\{
\begin{array}{c}
\partial _{t}\theta +\Lambda ^{\alpha }\theta +u\vec{\nabla}\theta =0 \\
u=\Lambda ^{\alpha -1}\mathcal{R}^{\perp }(\theta ) \\
\left. \theta \right\vert _{t=0}=\theta _{\ast }.%
\end{array}%
\right.  \tag{MQG}
\end{equation}%
Here, $\alpha \in ]0,1[$ is a \underline{fixed} real number, the unknown $%
\theta =\theta (t,x)$ is a real-valued function defined on $\mathbb{R}%
^{+}\times \mathbb{R}^{2},$ $\theta _{\ast }$ is a given initial data, $%
\mathcal{R}^{\perp }$ is the operator defined via Riesz transforms by;%
\begin{equation*}
\mathcal{R}^{\perp }(\theta )=\left( -\mathcal{R}_{2}\theta ,\mathcal{R}%
_{1}\theta \right) ,
\end{equation*}%
and $\Lambda ^{\gamma }$ is the non-local operator defined through the
Fourier transform by:%
\begin{equation*}
\widehat{\Lambda ^{\gamma }f}~(\xi )=\left\vert \xi \right\vert ^{\gamma }%
\hat{f}~(\xi ).
\end{equation*}%
The equation (MQG) was recently introduced in \cite{CIW} by P. Constantin,
G. Iyer and J. Wu as a modification of the 2D dissipative quasi-geostrophic
equation%
\begin{equation}
\left\{
\begin{array}{c}
\partial _{t}\theta +\Lambda ^{\alpha }\theta +u\vec{\nabla}\theta =0 \\
u=\mathcal{R}^{\perp }(\theta ) \\
\left. \theta \right\vert _{t=0}=\theta _{\ast }.%
\end{array}%
\right.  \tag{QG}
\end{equation}%
In \cite{CIW}, the autors proved that if the initial data $\theta _{\ast }$
belongs to $L^{2}(\mathbb{R}^{2})$ then the equation (MQG) has a global
solution $\theta \in C^{\infty }\left( ]0,+\infty \lbrack \times \mathbb{R}%
^{2}\right) .$

Very recently, C. Miao and L.\ Xue in \cite{MX} have proved that
for any initial
data $\theta _{\ast }\in H^{m}(\mathbb{R}^{2}),$ with $m\in \mathbb{N}$ and $%
m>2,$ there exists a unique global solution $\theta \in C(\mathbb{R}%
^{+},H^{m}(\mathbb{R}^{2}))$ to the equation (MQG). Moreover,the solution $%
\theta $ satisfies the following regularity property:%
\begin{equation*}
\forall \gamma \geq 0,~t^{\gamma }\theta \in L^{\infty }(\mathbb{R}%
^{+},H^{m+\gamma \sigma }(\mathbb{R}^{2})).
\end{equation*}%
One of the main properties of the equation (MQG) is the following scaling
invariance property: If $\theta $ is a solution to (MQG) with initial data $%
\theta _{\ast }$ then, for all $\lambda >0,$ the function $\theta _{\lambda
}\equiv \theta (\lambda ^{\alpha }t,\lambda x)$ is a solution to (MQG) with
initial data $\theta _{\ast }^{\lambda }\equiv \theta _{\ast }(\lambda x).$
This leads us to introduce the notion of critical space: a functional space $%
X$ is called a critical space for the equation (MQG) if for all $f\in X$ and
$\lambda >0,$ we have%
\begin{equation*}
\left\Vert f(\lambda .)\right\Vert _{X}=\left\Vert f\right\Vert _{X}.
\end{equation*}%
For instance, the homogenous Sobolev's space $\dot{H}^{\sigma }(\mathbb{R}%
^{2})$ is a critical space if and only if $\sigma =1.$ Therefore,
following the classical approach of Fujita-Kato \cite{FK}, it is
natural to ask if the equation (MQG) is well-posed if the initial
data $\theta _{\ast }$ belongs to the critical space
$\dot{H}^{1}(\mathbb{R}^{2}).$ In this paper, we give a slightly
weaker result. In fact, we prove the global well-posedness of the
smooth solution to the equation (MQG) when the initial data is in
the inhomogeneous Sobolev's space $H^{1}(\mathbb{R}^{2}).$
Precisely, our result states as follows

\begin{theorem}
\label{th1}Let $\theta _{\ast }\in H^{\sigma }(\mathbb{R}^{2})$ with $\sigma
\geq 1.$ Then there is a unique solution $\theta $ in
\begin{equation*}
C^{1}\left( \mathbb{R}^{+},L^{2}\left( \mathbb{R}^{2}\right) \right) \cap
C\left( \mathbb{R}^{+},H^{\sigma }\left( \mathbb{R}^{2}\right) \right) \cap
L_{loc}^{2}\left( \mathbb{R}^{+},H^{\sigma +\frac{\alpha }{2}}\left( \mathbb{%
R}^{2}\right) \right)
\end{equation*}%
to the equation (MQG). Moreover, for all $\sigma ^{\prime }\geq \sigma $ we
have%
\begin{equation*}
\theta \in C^{\infty }\left( \mathbb{R}_{\ast }^{+},H^{\sigma ^{\prime
}}\left( \mathbb{R}^{2}\right) \right) .
\end{equation*}%
In particular, $\theta \in C^{\infty }\left( \mathbb{R}_{\ast }^{+}\times
\mathbb{R}^{2}\right) .$
\end{theorem}

The proof of this theorem relies essentially on the following two
propositions. The first one is a local well-posedness result.

\begin{proposition}
\label{pr1}Let $\theta _{\ast }\in H^{\sigma }(\mathbb{R}^{2})$ with $\sigma
\geq 1.$ Then the equation (MQG) has a unique maximal solution $\theta $
belonging to the space
\begin{equation*}
C^{1}\left( [0,T^{\ast }[,L^{2}\left( \mathbb{R}^{2}\right) \right) \cap
C\left( [0,T^{\ast }[,H^{\sigma }\left( \mathbb{R}^{2}\right) \right) \cap
L_{loc}^{2}\left( [0,T^{\ast }[,H^{\sigma +\frac{\alpha }{2}}\left( \mathbb{R%
}^{2}\right) \right) .
\end{equation*}%
Moreover, the time $T^{\ast }$ is bounded from below by%
\begin{equation*}
\sup \{T>0:~\mathcal{K}(\theta _{\ast },T)\geq \varepsilon _{\sigma }\}
\end{equation*}%
where $\varepsilon _{\sigma }>0$ is a constant depending only on $\sigma ,$
and
\begin{equation*}
\mathcal{K}(\theta _{\ast },T)=\left\Vert \left( \left[ \frac{1-e^{-2\nu
2^{\alpha q}T}}{2\nu }\right] ^{1/2}2^{\sigma q}\left\Vert \Delta _{q}\theta
_{\ast }\right\Vert _{2}\right) _{q}\right\Vert _{l^{2}(\mathbb{Z})}
\end{equation*}%
where $\nu >0$ is an absolute constant and $\left( \Delta
_{q}\right) _{q}$ denotes the family of the Littlewood-Paley
operators (for the definition, see the section 2).
\end{proposition}

The second result concerns the propagation of the initial regularity.

\begin{proposition}
\label{pr2}Let $\theta _{\ast }\in H^{1}(\mathbb{R}^{2})$ and let $\theta $
be a solution to the equation (MQG) belonging to the space%
\begin{equation*}
C^{1}\left( [0,T],L^{2}\left( \mathbb{R}^{2}\right) \right) \cap C\left(
[0,T],H^{1}\left( \mathbb{R}^{2}\right) \right) \cap L^{2}\left( [0,T],H^{1+%
\frac{\alpha }{2}}\left( \mathbb{R}^{2}\right) \right) .
\end{equation*}%
If there exists $t_{0}\in \lbrack 0,T[$ and $\sigma \geq 1$ such that $%
\theta (t_{0})$ belongs to $H^{\sigma }(\mathbb{R}^{2}),$ then the solution $%
\theta $ belongs to the space%
\begin{equation*}
C\left( [t_{0},T],H^{\sigma }\left( \mathbb{R}^{2}\right) \right) \cap
L^{2}\left( [t_{0},T],H^{\sigma +\frac{\alpha }{2}}\left( \mathbb{R}%
^{2}\right) \right) .
\end{equation*}
\end{proposition}

\section{Notations and preliminaries}

In this preparatory section, we shall introduce some functionals spaces and
prove some elementary lemmas that will be used in the proof of Theorem \ref%
{th1}.

\subsection{Notations}

\begin{enumerate}
\item Throughout this paper, we will denotes various constants by
$C$. In particular, $C=C_{\ast ,\ast ,...}$ denotes constants
depending only on the quantities appearing in the index.

\item Let $A$ and $B$ be two reals functions. The notation $A\lesssim B$
means that there exists a constant $c>0$ such that $A\leq cB.$ We write $%
A\simeq B$ if $A\lesssim B$ and $B\lesssim A.$

\item For $p\in \lbrack 1,\infty ],$ we denote by $L^{p}$ the Lebesgue space
$L^{p}(\mathbb{R}^{2})$ endowed with the usual norm $\left\Vert .\right\Vert
_{p}.$

\item Let $T>0,~r\in \lbrack 1,\infty ]$ and $X$ be a Banach space. We
frequently denote the mixed space $L^{r}([0,T],X)$ by $L_{T}^{r}X.$

\item If $P$ and $Q$ are two operators, we denote by $[P,Q]$ the commutator
operator defined by%
\begin{equation*}
\lbrack P,Q]=PQ-QP.
\end{equation*}

\item $\langle .,.\rangle $ denotes the usual inner product in the
Hilbert space $L^{2}(\mathbb{R}^{2}).$

\item For $k\in \mathbb{N},$ we denote by $C_{B}^{k}(\mathbb{R}^{2})$ the
space of real-valued functions $f\in C^{k}(\mathbb{R}^{2})$ such that%
\begin{equation*}
\left\Vert f\right\Vert _{C_{B}^{k}(\mathbb{R}^{2})}\equiv \sup_{\beta \in
\mathbb{N}^{2},\left\vert \beta \right\vert \leq k}\left\Vert D^{\beta
}f\right\Vert _{\infty }<\infty .
\end{equation*}%
We set $C_{B}^{\infty }(\mathbb{R}^{2})=\cap _{k\in \mathbb{N}}C_{B}^{k}(%
\mathbb{R}^{2}).$
\end{enumerate}

\subsection{Sobolev's spaces and Chemin-Lerner's spaces}

We first recall the definition of the nonhomogeneous and homogenous Soblev's
spaces.

\begin{definition}
Let $s$ in $\mathbb{R}.$

\begin{enumerate}
\item The space $H^{s}(\mathbb{R}^{2})=H^{s}$ consists of all distributions $%
f\in S^{\prime }(\mathbb{R}^{2})$ such that
\begin{equation*}
\left\Vert f\right\Vert _{H^{s}}\equiv \left( \int_{\mathbb{R}^{2}}\left(
1+\left\vert \xi \right\vert ^{2}\right) ^{s}\left\vert \hat{f}(\xi
)\right\vert ^{2}d\xi \right) ^{1/2}<\infty .
\end{equation*}

\item The space $\dot{H}^{s}(\mathbb{R}^{2})=\dot{H}^{s}$ is the set of $%
f\in S^{\prime }(\mathbb{R}^{2})/\mathcal{P}(\mathbb{R}^{2})$ satisfying%
\begin{equation*}
\left\Vert f\right\Vert _{\dot{H}^{s}}\equiv \left( \int_{\mathbb{R}%
^{2}}\left\vert \xi \right\vert ^{2s}\left\vert \hat{f}(\xi )\right\vert
^{2}d\xi \right) ^{1/2}<\infty
\end{equation*}%
where $\mathcal{P}(\mathbb{R}^{2})$ is the space of polynomials functions
defined on $\mathbb{R}^{2}.$
\end{enumerate}
\end{definition}

\begin{remark}
Let $k\in \mathbb{N},~s\in \mathbb{R}$ and $\sigma \in \mathbb{R}^{+}.$ By
using the Plancherel formula, one can easily verify:%
\begin{eqnarray}
\left\Vert f\right\Vert _{\dot{H}^{s}} &\simeq &\left\Vert \Lambda
^{s}f\right\Vert _{2}  \label{eq1} \\
\left\Vert f\right\Vert _{H^{\sigma }} &\simeq &\left\Vert f\right\Vert
_{2}+\left\Vert \Lambda ^{\sigma }f\right\Vert _{2}  \label{eq2} \\
\left\Vert f\right\Vert _{H^{k}} &\simeq &\sum_{\beta \in \mathbb{N}%
^{2},\left\vert \beta \right\vert \leq k}\left\Vert D^{\beta }f\right\Vert
_{2}.  \label{eq3}
\end{eqnarray}
\end{remark}

\begin{notation}
Let $p\in \mathbb{N}\cup \{+\infty \}$ and $I\subset \mathbb{R}$. In the
sequel, we will use frequently the followings notations:%
\begin{eqnarray*}
H^{\infty } &\equiv &\cap _{k\in \mathbb{N}}H^{k} \\
C^{p}(I,H^{\infty }) &\equiv &\cap _{k\in \mathbb{N}}C^{p}(I,H^{k}).
\end{eqnarray*}
\end{notation}

Next, we recall the notion of the Littelwood-Paley decomposition that will
allow us to define the Chemin-Lerner spaces. Let $\phi \in S(\mathbb{R}^{2})$
such that $\hat{\phi}\equiv 1$ on the unit ball $B(0,1)$ of $\mathbb{R}^{2}$
and $\hat{\phi}\equiv 0$ outside the ball $B(0,2).$ For $q$ in $\mathbb{Z}$,
we denote by $S_{q}$ and $\Delta _{q}$ the operators defined by%
\begin{eqnarray*}
S_{q}f &=&\phi _{q}\ast f \\
\Delta _{q}f &=&S_{q+1}f-S_{q}f
\end{eqnarray*}%
where the star $\ast $ denotes the convolution on $\mathbb{R}^{2}$ and $\phi
_{q}=2^{2q}\phi (2^{q}.).$ For any $f$ in $S^{\prime }(\mathbb{R}^{2}),$ the
identity, called the Littelwood-Paley decomposition of $f$,%
\begin{equation}
f=\sum_{q\in \mathbb{Z}}\Delta _{q}f  \label{id}
\end{equation}%
holds in $S^{\prime }(\mathbb{R}^{2})/\mathcal{P}(\mathbb{R}^{2}).$
Moreover, if $f\in L^{p}$ with $1\leq p<\infty $ then the equality (\ref{id}%
) holds in $S^{\prime }(\mathbb{R}^{2}).$

It is well-known that the Littlewood-Paley decomposition provides
an
equivalent definition to the Sobolev spaces. Namely, we have%
\begin{eqnarray}
\left\Vert f\right\Vert _{H^{\sigma }} &\simeq &\left\Vert S_{0}f\right\Vert
_{2}+\left( \sum_{q\in \mathbb{N}}2^{2\sigma q}\left\Vert \Delta
_{q}f\right\Vert _{2}^{2}\right) ^{1/2}  \label{sob1} \\
\left\Vert f\right\Vert _{\dot{H}^{\sigma }} &\simeq &\left( \sum_{q\in
\mathbb{Z}}2^{2\sigma q}\left\Vert \Delta _{q}f\right\Vert _{2}^{2}\right)
^{1/2}  \label{sobolev2}
\end{eqnarray}

Now we recall the definition of Chemin-Lerner's spaces \cite{CL}

\begin{definition}
Let $T>0,~r\in \lbrack 1,\infty ]$ and $s\in \mathbb{R}.$

\begin{enumerate}
\item The space
$\tilde{L}^{r}([0,T],\dot{H}^{s}(\mathbb{R}^{2})),$ abbreviated by
$\tilde{L}_{T}^{r}\dot{H}^{s},$ is the set of all tempered
distribution $v$ satisfying%
\begin{equation*}
\left\Vert v\right\Vert _{\tilde{L}_{T}^{r}\dot{H}^{s}}\equiv \left(
\sum_{q\in \mathbb{Z}}2^{2sq}\left\Vert \Delta _{q}v\right\Vert
_{L_{T}^{r}L^{2}}^{2}\right) ^{1/2}<\infty .
\end{equation*}

\item If $s\geq 0,$ we denote by $\tilde{L}^{r}([0,T],H^{s}(\mathbb{R}%
^{2})), $ abbreviated by $\tilde{L}_{T}^{r}H^{s},$ the space $%
L_{T}^{r}L^{2}\cap \tilde{L}_{T}^{r}\dot{H}^{s}$ endowed with the norm%
\begin{equation*}
\left\Vert v\right\Vert _{\tilde{L}_{T}^{r}H^{s}}\equiv \left\Vert
v\right\Vert _{L_{T}^{r}L^{2}}+\left\Vert v\right\Vert _{\tilde{L}_{T}^{r}%
\dot{H}^{s}}.
\end{equation*}
\end{enumerate}
\end{definition}

Using the estimate (\ref{sobolev2}) and the Minkowiski inequality one can
can easily obtain the following estimates%
\begin{eqnarray}
\left\Vert v\right\Vert _{\tilde{L}_{T}^{r}\dot{H}^{s}} &\lesssim
&\left\Vert v\right\Vert _{L_{T}^{r}\dot{H}^{s}}\text{ if }r\leq 2
\label{min1} \\
\left\Vert v\right\Vert _{L_{T}^{r}\dot{H}^{s}} &\lesssim &\left\Vert
v\right\Vert _{\tilde{L}_{T}^{r}\dot{H}^{s}}\text{ if }r\geq 2  \label{min2}
\\
\left\Vert v\right\Vert _{L_{T}^{r}\dot{H}^{s}} &\simeq &\left\Vert
v\right\Vert _{\tilde{L}_{T}^{r}\dot{H}^{s}}\text{ if }r=2.  \label{min3}
\end{eqnarray}

In the sequel we will often use the following notation.

\begin{notation}
Let $\sigma \geq 0$ and $T>0.$ We set%
\begin{eqnarray}
\mathbf{X}_{T}^{\sigma } &\equiv &\tilde{L}_{T}^{\infty }H^{\sigma }\cap
L_{T}^{2}H^{\sigma +\frac{\alpha }{2}}  \label{x} \\
\mathbf{Z}_{T}^{\sigma } &\equiv &C^{1}([0,T],L^{2})\cap C([0,T],H^{\sigma
})\cap L_{T}^{2}H^{\sigma +\frac{\alpha }{2}}.  \label{z}
\end{eqnarray}
\end{notation}

\subsection{Preliminaries results}

In this subsection, we state and prove some elementary and useful
lemmas.

The first lemma is a particular case of the well-known Bernstein inequality
\cite{L}.

\begin{lemma}
\label{lem1}Let $\gamma \in \mathbb{R},~\beta \in \mathbb{N}^{2}$ and $1\leq
p\leq r\leq \infty .$ The followings assertions hold true:

\begin{enumerate}
\item There exists a constant $\nu =\nu (\gamma )>0$ such that for all $f\in
S^{\prime }(\mathbb{R}^{2})$ and $q\in \mathbb{Z}$ we have%
\begin{equation}
\left\Vert \Lambda ^{\gamma }\Delta _{q}f\right\Vert _{p}\geq \nu 2^{\gamma
q}\left\Vert \Delta _{q}f\right\Vert _{p}.  \label{bern1}
\end{equation}

\item There exists a constant $C=C(\gamma ,\beta ,p,r)>0$ such that for all $%
f\in S^{\prime }(\mathbb{R}^{2})$ and $q\in \mathbb{Z}$ we have%
\begin{equation}
\left\Vert \Lambda ^{\gamma }D^{\beta }\Delta _{q}f\right\Vert _{r}\leq
C2^{q(\gamma +\left\vert \beta \right\vert )+2q\left( \frac{1}{p}-\frac{1}{r}%
\right) }\left\Vert \Delta _{q}f\right\Vert _{p}.  \label{bern2}
\end{equation}
\end{enumerate}
\end{lemma}

\bigskip

The next lemma will be repeatedly used in this paper

\begin{lemma}
\label{lem2}Let $\gamma \in ]0,1[,~k\in \mathbb{N}$ and $m\in \mathbb{N}%
\backslash \{0,1,2\}.$ We have the following assertions:

\begin{enumerate}
\item For all $f\in H^{k}(\mathbb{R}^{2})$ and $g\in C_{B}^{k}(\mathbb{R}%
^{2}),$%
\begin{equation}
\left\Vert fg\right\Vert _{H^{k}}\lesssim \left\Vert g\right\Vert
_{C_{B}^{k}(\mathbb{R}^{2})}\left\Vert f\right\Vert _{H^{k}}.  \label{eg1}
\end{equation}

\item For all $f\in L^{2}(\mathbb{R}^{2}),$%
\begin{equation}
\left\Vert \Lambda ^{\gamma -1}\mathcal{R}^{\perp }(f)\right\Vert _{\frac{2}{%
\gamma }}\lesssim \left\Vert f\right\Vert _{2}.  \label{eg2}
\end{equation}

\item For all $f\in H^{m}(\mathbb{R}^{2}),$%
\begin{equation}
\left\Vert \Lambda ^{\gamma -1}\mathcal{R}^{\perp }(f)\right\Vert
_{C_{B}^{m-2}(\mathbb{R}^{2})}\lesssim \left\Vert f\right\Vert _{H^{m}}.
\label{eg3}
\end{equation}%
In particular, if $f\in H^{\infty }(\mathbb{R}^{2})$ then $\Lambda ^{\gamma }%
\mathcal{R}^{\perp }(f)\in C_{B}^{\infty }(\mathbb{R}^{2}).$

\item For all $f\in C_{B}^{k}(\mathbb{R}^{2})\cap L^{2}(\mathbb{R}^{2}),$%
\begin{equation}
\left\Vert \Lambda ^{\gamma -1}\mathcal{R}^{\perp }(f)\right\Vert
_{C_{B}^{k}(\mathbb{R}^{2})}\lesssim \left( \left\Vert f\right\Vert
_{2}+\left\Vert f\right\Vert _{C_{B}^{k}(\mathbb{R}^{2})}\right) .
\label{eg4}
\end{equation}
\end{enumerate}
\end{lemma}

\begin{proof}
(1) obvious. (2) is a direct consequence of the Sobolev embedding $\dot{H}%
^{1-\gamma }\left( \mathbb{R}^{2}\right) \hookrightarrow L^{\frac{2}{\gamma }%
}(\mathbb{R}^{2})$ and the continuity of the Riesz transforms on $L^{2}(%
\mathbb{R}^{2}).$ Let us prove the assertion (3). Let $f\in H^{m}(\mathbb{R}%
^{2})$ and set $g\equiv \Lambda ^{\gamma -1}\mathcal{R}^{\perp }(f).$
According to (\ref{eg2}), $g\in L^{\frac{2}{\gamma }}(\mathbb{R}^{2}).$
Consequentlty,
\begin{equation*}
g=\sum_{q\in Z}\Delta _{q}g\text{ in }S^{\prime }(\mathbb{R}^{2}).
\end{equation*}%
Hence, for all $\beta \in \mathbb{N}^{2}$ with $\left\vert \beta \right\vert
\leq m-2,$ we have, thanks to (\ref{bern2}),
\begin{eqnarray*}
\left\Vert D^{\beta }g\right\Vert _{\infty } &\lesssim &\sum_{q\leq
0}2^{q(\gamma +\left\vert \beta \right\vert )}\left\Vert \Delta
_{q}f\right\Vert _{2}+\sum_{q>0}2^{q(\gamma -1)}\left\Vert \Delta
_{q}D^{\beta }f\right\Vert _{\infty } \\
&\lesssim &\left\Vert f\right\Vert _{2}+\left\Vert D^{\beta }f\right\Vert
_{\infty } \\
&\lesssim &\left\Vert f\right\Vert _{H^{m}}.
\end{eqnarray*}%
In the last inequality, we have used the embedding $H^{m}(\mathbb{R}%
^{2})\hookrightarrow C_{B}^{m-2}(\mathbb{R}^{2}).$ The proof of the
assertion (4) is similar to that of the assertion (3) and thus it is omitted.
\end{proof}

\bigskip

The next lemma is classical. For the proof see for instance \cite{RS} or
\cite{L}

\begin{lemma}
\label{lem3}Let $\sigma _{1}$ and $\sigma _{2}$ be two real numbers such
that $\sigma _{1}<1,\sigma _{2}<1$ and $\sigma _{1}+\sigma _{2}>0.$ Then
there exists a constant $C=C_{\sigma _{1},\sigma _{2}}\geq 0$ such that for
all $f\in \dot{H}^{\sigma _{1}}\left( \mathbb{R}^{2}\right) $ and $g\in \dot{%
H}^{\sigma _{2}}\left( \mathbb{R}^{2}\right) $ we have%
\begin{equation}
\left\Vert fg\right\Vert _{\dot{H}^{\sigma }}\leq C\left\Vert f\right\Vert _{%
\dot{H}^{\sigma _{1}}}\left\Vert g\right\Vert _{\dot{H}^{\sigma _{2}}}
\label{mul}
\end{equation}%
where $\sigma =\sigma _{1}+\sigma _{2}-1.$
\end{lemma}

\bigskip

The following result on differential inequalities will be useful.

\begin{lemma}
\label{lem4}Let $T>0$ and $h,f:[0,T]\rightarrow \mathbb{R}^{+}$ two
continuous functions. Assume that $h^{2}\in C^{1}\left( [0,T]\right) $ and
there exists $c\in \mathbb{R}$ such that for all $t$ in $[0,T]$%
\begin{equation*}
\left( h^{2}\right) ^{\prime }(t)+c~h^{2}(t)\leq f(t)h(t).
\end{equation*}%
Then%
\begin{equation}
\forall t\in \lbrack 0,T],~h(t)\leq e^{-\frac{c}{2}t}h(0)+\frac{1}{2}%
\int_{0}^{t}e^{-c(t-s)}f(s)ds.  \label{gr}
\end{equation}
\end{lemma}

\begin{proof}
For $\varepsilon >0,$ we define the function $h_{\varepsilon }=\sqrt{%
\varepsilon +h^{2}}.$ Clearly, $h_{\varepsilon }$ is $C^{1}$ on $[0,T]$ and
satisfies the differential inequality%
\begin{eqnarray*}
h_{\varepsilon }^{\prime }+\frac{c}{2}h_{\varepsilon } &\leq &\frac{%
c\varepsilon }{2\sqrt{\varepsilon +h^{2}}}+\frac{fh}{2\sqrt{\varepsilon
+h^{2}}} \\
&\leq &\frac{c\sqrt{\varepsilon }}{2}+\frac{1}{2}f.
\end{eqnarray*}%
Then for all $t$ in $[0,T]$ we have%
\begin{equation*}
\ h_{\varepsilon }(t)\leq e^{-\frac{c}{2}t}h_{\varepsilon }(0)+\frac{c}{2}%
\sqrt{\varepsilon }\int_{0}^{t}e^{-c(t-s)}ds+\frac{1}{2}%
\int_{0}^{t}e^{-c(t-s)}f(s)ds.
\end{equation*}%
Hence, we get the desired estimate (\ref{gr}) by sending $\varepsilon
\rightarrow 0.$
\end{proof}

\bigskip

We now state and prove a version of the well-known maximal principle.

\begin{lemma}
\label{lem5}Let $T>0,~\theta \in C^{1}([0,T],L^{2})\cap C([0,T],H^{1})$ and $%
f,v\in C([0,T],L^{2})$ such that
\begin{equation*}
\partial _{t}\theta +\Lambda ^{\alpha }\theta +u\vec{\nabla}\theta =f
\end{equation*}%
where $u=\Lambda ^{\alpha -1}\mathcal{R}^{\perp }(v).$ Then, for all $t$ in $%
[0,T],$ we have%
\begin{equation}
\frac{1}{2}\frac{d}{dt}\left\Vert \theta (t)\right\Vert _{2}^{2}+\left\Vert
\Lambda ^{\alpha /2}\theta (t)\right\Vert _{2}^{2}=\langle f(t),\theta
(t)\rangle  \label{max1}
\end{equation}%
and
\begin{equation}
\left\Vert \theta (t)\right\Vert _{2}\leq \left\Vert \theta (0)\right\Vert
_{2}+\int_{0}^{t}\left\Vert f(\tau )\right\Vert _{2}d\tau .  \label{max2}
\end{equation}
\end{lemma}

\begin{proof}
Let $t$ in $[0,T].$ We have%
\begin{eqnarray*}
\frac{1}{2}\frac{d}{dt}\left\Vert \theta (t)\right\Vert _{2}^{2} &=&\langle
\partial _{t}\theta (t),\theta (t)\rangle \\
&=&-\underset{=I(t)}{\underbrace{\langle \Lambda ^{\alpha }\theta (t),\theta
(t)\rangle }}-\underset{=J(t)}{\underbrace{\langle u(t)\vec{\nabla}\theta
(t),\theta (t)\rangle }}+\langle f(t),\theta (t)\rangle .
\end{eqnarray*}%
By Plancherel's formula,%
\begin{equation*}
I(t)=\left\Vert \Lambda ^{\alpha /2}\theta (t)\right\Vert _{2}^{2}.
\end{equation*}%
On the other hand, by virtue of the density, there exists two sequence $%
(v_{n})_{n}$ and $(\theta _{n})_{n}$ in $C_{c}^{\infty }(\mathbb{R}^{2})$
such that
\begin{equation*}
v_{n}\rightarrow v(t)\text{ in }L^{2}\text{ and }\theta _{n}\rightarrow
\theta (t)\text{ in }H^{1}.
\end{equation*}%
Consequently, the estimate (\ref{eg2}) and the Sobolev embedding $%
H^{1}\hookrightarrow L^{\frac{2}{1-\alpha }}$ imply that%
\begin{equation*}
u_{n}\equiv \Lambda ^{\alpha -1}\mathcal{R}^{\perp }(v_{n})\rightarrow u%
\text{ in }L^{\frac{2}{\alpha }}\text{ and }\theta _{n}\rightarrow \theta (t)%
\text{ in }L^{\frac{2}{1-\alpha }}.
\end{equation*}%
Hence, Holder's inequality yields%
\begin{equation*}
J_{n}\equiv \langle u_{n}\vec{\nabla}\theta _{n},\theta _{n}\rangle
\rightarrow J(t)\text{.}
\end{equation*}%
Now, a simple integration by parts using the fact that $u_{n}$ is
divergence-free gives%
\begin{equation*}
J_{n}=-J_{n}
\end{equation*}%
Thus $J_{n}=0$ and consequently $J(t)=0.$ This finishes the proof of (\ref%
{max1}). Finally, the estimate (\ref{max2}) is an immediate consequence of (%
\ref{max1}) and Lemma \ref{lem5}.
\end{proof}

\bigskip

The following lemma is simple, hence its proof is omitted.

\begin{lemma}
\label{lem6}Let $T>0$, $\sigma \in \mathbb{R}$ and $r\in \lbrack 1,\infty ].$
Assume that a sequence $(\theta _{n})_{n}$ converges to a function $\theta $
in the space $L^{\infty }([0,T],L^{2}(\mathbb{R}^{2}))$ and that
\begin{equation*}
\mathcal{M}\equiv \sup_{n}\left\Vert \theta _{n}\right\Vert _{\tilde{L}%
_{T}^{r}\dot{H}^{\sigma }}<\infty .
\end{equation*}%
Then $\theta \in \tilde{L}_{T}^{r}\dot{H}^{\sigma }$ and satisfies%
\begin{equation*}
\left\Vert \theta \right\Vert _{\tilde{L}_{T}^{r}\dot{H}^{\sigma }}\leq
\mathcal{M}.
\end{equation*}
\end{lemma}

\begin{lemma}
\label{lem7}Let $T>0$ and $\sigma \in \mathbb{R}$. If a function $\theta $
belongs to the spaces $C([0,T],L^{2}(\mathbb{R}^{2}))$ and $\tilde{L}%
^{\infty }\left( [0,T],\dot{H}^{\sigma }\right) $ then it belongs to the
space $C\left( [0,T],\dot{H}^{\sigma }\right) .$
\end{lemma}

\begin{proof}
For $N\in \mathbb{N}$ define $\theta _{N}=\sum_{\left\vert j\right\vert \leq
N}\Delta _{j}\theta .$ The sequence $\left( \theta _{N}\right) _{N}$ is in
the space $C\left( [0,T],\dot{H}^{\sigma }\right) ;$ in fact, for every $%
t,t^{\prime }\in \lbrack 0,T]$ we have%
\begin{eqnarray*}
\left\Vert \theta _{N}(t)-\theta _{N}(t^{\prime })\right\Vert _{\dot{H}%
^{\sigma }}^{2} &=&\sum_{\left\vert j\right\vert \leq N}2^{2\sigma
j}\left\Vert \Delta _{j}\left[ \theta (t)-\theta (t^{\prime })\right]
\right\Vert _{2}^{2} \\
&\leq &\left\Vert \theta (t)-\theta (t^{\prime })\right\Vert
_{2}\sum_{\left\vert j\right\vert \leq N}2^{2\sigma j}.
\end{eqnarray*}%
On the other hand, $\left( \theta _{N}\right) _{N}$ converges to $\theta $
in the space $L^{\infty }\left( [0,T],\dot{H}^{\sigma }\right) ;$ indeed,
for every $t$ in the interval $[0,T]$ we have%
\begin{equation*}
\left\Vert \theta (t)-\theta _{N}(t)\right\Vert _{\dot{H}^{\sigma }}^{2}\leq
4\sum_{\left\vert j\right\vert >N}2^{2\sigma j}\left\Vert \Delta _{j}\theta
\right\Vert _{L_{T}^{\infty }L^{2}}^{2}
\end{equation*}%
and $\sum_{\left\vert j\right\vert >N}2^{2\sigma j}\left\Vert \Delta
_{j}\theta \right\Vert _{L_{T}^{\infty }L^{2}}^{2}\rightarrow 0$ as $%
N\rightarrow \infty $ since $\theta \in \tilde{L}^{\infty }\left( [0,T],\dot{%
H}^{\sigma }\right) .$ Finally, since $C\left( [0,T],\dot{H}^{\sigma
}\right) $ is closed sub-space of $L^{\infty }\left( [0,T],\dot{H}^{\sigma
}\right) ,$ we conclude that $\theta $ is in $C\left( [0,T],\dot{H}^{\sigma
}\right) .$
\end{proof}

\bigskip

Now we state and prove the main result of this sub-section.

\begin{proposition}
\label{pr3}Let $\alpha \in ]0,1[,~T>0,~\theta _{\ast }\in H^{\infty }(%
\mathbb{R}^{2})$ and $v\in C([0,T],H^{\infty }(\mathbb{R}^{2})).$ Set $%
u=\Lambda ^{\alpha -1}\mathcal{R}^{\perp }(v).$ Then there is a unique
solution $\theta \in C^{1}([0,T],H^{\infty }(\mathbb{R}^{2}))$ to the linear
initial value problem%
\begin{equation}
\left\{
\begin{array}{c}
\partial _{t}\theta +\Lambda ^{\alpha }\theta +u\vec{\nabla}\theta =0 \\
\left. \theta \right\vert _{t=0}=\theta _{\ast }.%
\end{array}%
\right.  \tag{IVPL}
\end{equation}%
Moreover,%
\begin{equation}
\sup_{0\leq t\leq T}\left\Vert \theta (t)\right\Vert _{2}\leq \left\Vert
\theta _{\ast }\right\Vert _{2}  \label{r1}
\end{equation}%
and, for every $r\geq 2$ and $s\geq 0$ there exists a constant $\mathcal{C}%
_{r,s}>0$ depending only on $r$ and $s$ such that%
\begin{equation}
\left\Vert \theta \right\Vert _{\tilde{L}_{T}^{r}\dot{H}^{s+\frac{\alpha }{r}%
}}\leq \mathcal{K}_{r,s}(\theta _{\ast },T)+\mathcal{C}_{r,s}\left(
\left\Vert v\right\Vert _{L_{T}^{2}\dot{H}^{1+\frac{\alpha }{2}}}\left\Vert
\theta \right\Vert _{\tilde{L}_{T}^{r}\dot{H}^{s+\frac{\alpha }{r}%
}}+\left\Vert \theta \right\Vert _{L_{T}^{2}\dot{H}^{1+\frac{\alpha }{2}%
}}\left\Vert v\right\Vert _{\tilde{L}_{T}^{r}\dot{H}^{s+\frac{\alpha }{r}%
}}\right)  \label{r2}
\end{equation}%
where%
\begin{equation}
\mathcal{K}_{r,s}(\theta _{\ast },T)=\left\Vert \left( \frac{1-e^{-\nu
r2^{\alpha q}T}}{\nu r}\right) ^{1/r}2^{sq}\left\Vert \Delta _{q}\theta
_{\ast }\right\Vert _{2}\right\Vert _{l^{2}(\mathbb{Z})}\text{ if }r<\infty ,
\label{r3}
\end{equation}%
and%
\begin{equation}
\mathcal{K}_{\infty ,s}(\theta _{\ast },T)=\left\Vert 2^{sq}\left\Vert
\Delta _{q}\theta _{\ast }\right\Vert _{2}\right\Vert _{l^{2}(\mathbb{Z})}.
\label{inf}
\end{equation}%
$\nu $ is the real given by Lemma \ref{lem2}.
\end{proposition}

\bigskip

The proof of this proposition is based on the following commutator
estimate which can be easily proved by following the arguments
used for instance in \cite{C}, \cite{CMZ}, \cite{D1}, \cite{Mura}
or \cite{W} .

\begin{lemma}[Commutator estimate]
\label{lem8}Let $\alpha \in ]0,1[,~T>0$ and $f,g$ two functions in the space
$C([0,T],H^{\infty }(\mathbb{R}^{2})).$ Set $u=\Lambda ^{\alpha -1}\mathcal{R%
}^{\perp }(f).$ Then for all $r\geq 2$ and $s\geq 0$ there exists a constant
$C_{r,s}$ depending only on $r$ and $s$ and a sequence $\left( \varepsilon
_{q}\right) _{q}\in l^{2}(\mathbb{Z})$ with $\sum_{q}\varepsilon
_{q}^{2}\leq 1$ such that for all $q\in \mathbb{Z}$ we have%
\begin{equation}
\left\Vert \left[ \Delta _{q},u\right] \vec{\nabla}g\right\Vert
_{L_{T}^{\rho }L^{2}}\leq C_{r,s}2^{-q(s+\frac{\alpha }{r}-\frac{\alpha }{2}%
)}\varepsilon _{q}\left( \left\Vert f\right\Vert _{L_{T}^{2}\dot{H}^{1+\frac{%
\alpha }{2}}}\left\Vert g\right\Vert _{\tilde{L}_{T}^{r}\dot{H}^{s+\frac{%
\alpha }{r}}}+\left\Vert g\right\Vert _{L_{T}^{2}\dot{H}^{1+\frac{\alpha }{2}%
}}\left\Vert f\right\Vert _{\tilde{L}_{T}^{r}\dot{H}^{s+\frac{\alpha }{r}%
}}\right)  \label{com}
\end{equation}%
where $\frac{1}{\rho }=\frac{1}{2}+\frac{1}{r}.$
\end{lemma}

\begin{proof}[Proof of Proposition 3]
We first notice that thanks to Lemma \ref{lem2}, the function $u$ belongs to
the space $\cap _{k\in \mathbb{N}}C([0,T],C_{B}^{k}(\mathbb{R}^{2})).$
Hence, for all $m\in \mathbb{N},$%
\begin{equation}
\mathcal{M}_{m}(u)\equiv \sup_{\left\vert \beta \right\vert \leq
m}\left\Vert D^{\beta }u\right\Vert _{L^{\infty }([0,T]\times \mathbb{R}%
^{2})}<\infty .  \label{uni}
\end{equation}%
Let us prove the existence of the solution $\theta .$ To do so, we will make
use of the classical Fredireck's method. For $n\in \mathbb{N},$ we consider
the linear ODE%
\begin{equation}
\left\{
\begin{array}{c}
\partial _{t}\theta =F_{n}(t,\theta ) \\
\left. \theta \right\vert _{t=0}=J_{n}\theta _{\ast }%
\end{array}%
\right.  \tag{S$_{n}$}
\end{equation}%
where the operator $J_{n}$ is defined by%
\begin{equation*}
\widehat{J_{n}f}(\xi )=1_{B(0,n)}(\xi )~\hat{f}(\xi )
\end{equation*}%
and $F_{n}:[0,T]\times L^{2}(\mathbb{R}^{2})\rightarrow L^{2}(\mathbb{R}%
^{2}) $ is the function defined by%
\begin{equation*}
F_{n}\left( t,g\right) \equiv -\Lambda ^{\alpha }J_{n}g-J_{n}(u\vec{\nabla}%
J_{n}\theta ).
\end{equation*}%
Clearly, $F_{n}$ belongs to $C([0,T],\mathcal{L}\left( L^{2}(\mathbb{R}%
^{2})\right) ).$ Hence the Cauchy-Lipschitz Theorem ensures that
(S$_{n}$)
has a unique global solution $\theta _{n}$ belonging to the space $%
C^{1}([0,T],L^{2}(\mathbb{R}^{2})).$ Now, since $J_{n}^{2}=J_{n},$ then $%
J_{n}\theta _{n}$ is also a solution to (S$_{n}$). Therefore, $\theta _{n}={J%
}_{n}\theta _{n}$ which implies in particular that $\theta _{n}\in
C^{1}([0,T],H^{\infty }(\mathbb{R}^{2})).$ Now we will estimate the norm of $%
\theta _{n}$ in the spaces $L_{T}^{\infty }H^{m}$ where $m\in \mathbb{N}.$
Let $m\in \mathbb{N}$ and $\beta \in \mathbb{N}^{2}$ such that $\left\vert
\beta \right\vert \leq m.$ Clearly, the function $D^{\beta }\theta _{n}$
satisfies the equation%
\begin{equation*}
\partial _{t}D^{\beta }\theta _{n}+\Lambda ^{\alpha }D^{\beta }\theta _{n}+u%
\vec{\nabla}D^{\beta }\theta _{n}=J_{n}\left( \left[ u,D^{\beta }\right]
\vec{\nabla}\theta _{n}\right) +\tilde{J}_{n}\left( u\vec{\nabla}D^{\beta
}\theta _{n}\right)
\end{equation*}%
where $\tilde{J}_{n}$ is defined by%
\begin{equation}
\widehat{\tilde{J}_{n}f}(\xi )=\left( 1-1_{B(0,n)}(\xi )\right) ~\hat{f}(\xi
).  \label{a}
\end{equation}%
Hence, Lemma \ref{lem5} implies%
\begin{eqnarray*}
\frac{d}{2dt}\left\Vert D^{\beta }\theta _{n}(t)\right\Vert _{2}^{2} &\leq
&\langle J_{n}\left( \left[ u,D^{\beta }\right] \vec{\nabla}\theta
_{n}\right) ,D^{\beta }\theta _{n}\rangle +\langle \tilde{J}_{n}\left( u\vec{%
\nabla}D^{\beta }\theta _{n}\right) ,D^{\beta }\theta _{n}\rangle \\
&=&\langle \left[ u,D^{\beta }\right] \vec{\nabla}\theta _{n},D^{\beta
}\theta _{n}\rangle \\
&\leq &\left\Vert \left[ u,D^{\beta }\right] \vec{\nabla}\theta
_{n}\right\Vert _{2}\left\Vert D^{\beta }\theta _{n}(t)\right\Vert _{2}^{2}
\\
&\leq &C_{\beta }~\mathcal{M}_{m}(u)\left\Vert \theta _{n}(t)\right\Vert
_{H^{m}}^{2}
\end{eqnarray*}%
where in the second equation we have used the fact $J_{n}D^{\beta }\theta
_{n}=D^{\beta }\theta _{n}$ and $\tilde{J}_{n}D^{\beta }\theta _{n}=0.$ We
then deduce that there exists a constant $C_{m}>0$ depending only on $m,$
such that
\begin{equation*}
\frac{d}{dt}\left\Vert \theta _{n}(t)\right\Vert _{H^{m}}^{2}\leq C_{m}%
\mathcal{M}_{m}(u)\left\Vert \theta _{n}(t)\right\Vert _{H^{m}}^{2}
\end{equation*}%
Invoking Gronwall's inequality, we then infer that the sequence $(\theta
_{n})_{n}$ is bounded in the space $C([0,T],H^{m}(\mathbb{R}^{2}))$. That is,%
\begin{equation}
\forall m\in \mathbb{N},~\Gamma _{m}\equiv \sup_{n}\left\Vert \theta
_{n}\right\Vert _{L_{T}^{\infty }H^{m}}<\infty .  \label{uniforme}
\end{equation}%
Next we will show that $(\theta _{n})_{n}$ is of Cauchy in the space $%
C([0,T],L^{2}(\mathbb{R}^{2})).$ Let $p,q\in \mathbb{N}$ such that $p\leq q.$
The function $\omega =\theta _{q}-\theta _{p}$ satisfies the equation%
\begin{equation*}
\partial _{t}\omega +\Lambda ^{\alpha }\omega +u\vec{\nabla}\omega =\left(
J_{q}-J_{p}\right) \left( u\vec{\nabla}\theta _{p}\right) +\tilde{J}%
_{q}\left( u\vec{\nabla}\omega \right)
\end{equation*}%
where $\tilde{J}_{q}$ is defined by (\ref{a}). Using Lemma \ref{lem5} and the fact that $\tilde{J}%
_{q}\omega =0,$ we easily get the following estimates%
\begin{eqnarray*}
\frac{d\left\Vert \omega (t)\right\Vert _{2}^{2}}{2dt} &\leq &\left\Vert
\left( J_{q}-J_{p}\right) \left( u\vec{\nabla}\theta _{p}\right) \right\Vert
_{2}\left\Vert \omega (t)\right\Vert _{2} \\
&\lesssim &\frac{1}{p}\left\Vert \left( u\vec{\nabla}\theta _{p}\right)
\right\Vert _{H^{1}}\left\Vert \omega (t)\right\Vert _{2} \\
&\lesssim &\frac{1}{p}\mathcal{M}_{2}(u)\left\Vert \theta _{p}\right\Vert
_{H^{2}}\left\Vert \omega (t)\right\Vert _{2} \\
&\lesssim &\frac{1}{p}\left\Vert \omega (t)\right\Vert _{2}.
\end{eqnarray*}%
In the last inequality, we used (\ref{uni})-(\ref{uniforme}). Hence, Lemma %
\ref{lem4} implies that there exists a constant $C$ independent on $p$ and $%
q $ such that%
\begin{equation*}
\sup_{0\leq t\leq T}\left\Vert \theta _{q}(t)-\theta _{p}(t)\right\Vert
_{2}\leq \left\Vert \left( J_{q}-J_{p}\right) \theta _{\ast }\right\Vert
_{2}+\frac{C}{p},
\end{equation*}%
which leads the required result. Let $\theta $ be the limit of the sequence $%
\left( \theta _{n}\right) _{n}$ in the space $C([0,T],L^{2}(\mathbb{R}%
^{2})). $ Now thanks to the interpolation inequality%
\begin{equation*}
\left\Vert f\right\Vert _{H^{m}}\leq C_{m}\sqrt{\left\Vert f\right\Vert
_{2}\left\Vert f\right\Vert _{H^{2m}}}
\end{equation*}%
and the uniform boundness (\ref{uniforme}), we infer that the sequence $%
\left( \theta _{n}\right) _{n}$ converges to $\theta $ in the space $%
C([0,T],H^{m}(\mathbb{R}^{2}))$ for all $m$ in $\mathbb{N}.$ Moreover, since
\begin{equation*}
\partial _{t}\theta _{n}=-\Lambda ^{\alpha }\theta _{n}-J_{n}\left( u\vec{%
\nabla}\theta _{n}\right)
\end{equation*}%
then by using the first assertion of Lemma \ref{lem1} and the continuity of
the operators $\Lambda ^{\alpha }$ and $\vec{\nabla}$ from $H^{m}$ into $%
H^{m-1},$ we easily deduce that, for all for all $m$ in $\mathbb{N},$ the
sequence $\left( \partial _{t}\theta _{n}\right) _{n}$ converges to $%
-\Lambda ^{\alpha }\theta -u\vec{\nabla}\theta $ in the space $C([0,T],H^{m}(%
\mathbb{R}^{2})).$ We therefore conclude that $\theta $ belongs to $%
C^{1}([0,T],H^{\infty }(\mathbb{R}^{2}))$ and it is a solution to (IVPL).The
uniqueness can be easily proved, indeed if $\theta ^{\prime }$ is another
solution to (IVPL) then the difference function $\delta =\theta -\theta
^{\prime }$ satisfies%
\begin{equation*}
\left\{
\begin{array}{c}
\partial _{t}\delta +\Lambda ^{\alpha }\delta +u\vec{\nabla}\delta =0 \\
\left. \delta \right\vert _{t=0}=0,%
\end{array}%
\right.
\end{equation*}%
which en virtue of Lemma (\ref{lem5}) implies $\delta =0$ and consequently, $%
\theta =\theta ^{\prime }.$

Once again the estimate (\ref{r1}) is a consequence of Lemma \ref{lem5}.
Finally, let us prove the estimate (\ref{r2}) in the case $r\in \lbrack
2,\infty \lbrack ,$ the proof in the case $r=\infty $ is similar and even
more simpler.

Apply the operator $\Delta _{q}$ ($q\in \mathbb{Z})$ to the first equation
of (IVPL), we get%
\begin{equation*}
\partial _{t}\Delta _{q}\theta +\Lambda ^{\alpha }\Delta _{q}\theta +u\vec{%
\nabla}\Delta _{q}\theta =\left[ u,\Delta _{q}\right] \vec{\nabla}\theta .
\end{equation*}%
Therefore, Lemma \ref{lem5}\ implies%
\begin{equation*}
\frac{d}{2dt}\left\Vert \Delta _{q}\theta (t)\right\Vert _{2}^{2}+\left\Vert
\Lambda ^{\alpha /2}\Delta _{q}\theta (t)\right\Vert _{2}^{2}\leq
F_{q}(t)\left\Vert \Delta _{q}\theta (t)\right\Vert _{2}
\end{equation*}%
where
\begin{equation*}
F_{q}(t)=\left\Vert \left[ u,\Delta _{q}\right] \vec{\nabla}\theta
\right\Vert _{2}.
\end{equation*}%
Thanks to the first assertion of Lemma \ref{lem1}, we deduce that there
exists a pure constant $\nu >0$ such that%
\begin{equation*}
\frac{d}{2dt}\left\Vert \Delta _{q}\theta (t)\right\Vert _{2}^{2}+\nu
2^{\alpha q}\left\Vert \Delta _{q}\theta (t)\right\Vert _{2}^{2}\leq
F_{q}(t)\left\Vert \Delta _{q}\theta (t)\right\Vert _{2}.
\end{equation*}%
Invoking Lemma \ref{lem4}, we obtain%
\begin{equation*}
\left\Vert \Delta _{q}\theta (t)\right\Vert _{2}\leq e^{-\nu 2^{\alpha
q}t}\left\Vert \Delta _{q}\theta _{\ast }\right\Vert
_{2}+\int_{0}^{t}e^{-\nu 2^{\alpha q}(t-s)}F_{q}(s)ds.
\end{equation*}%
Let $r\in \lbrack 2,\infty \lbrack $ and set $\rho =(\frac{1}{2}+\frac{1}{r}%
)^{-1}.$ Using the Young inequality, we deduce from the above inequality
that
\begin{eqnarray*}
\left\Vert \Delta _{q}\theta \right\Vert _{L_{T}^{r}L^{2}} &\leq &\left(
\frac{1-e^{-\nu 2^{\alpha q}T}}{\nu r}\right) ^{1/r}2^{-\frac{\alpha }{r}%
q}\left\Vert \Delta _{q}\theta _{\ast }\right\Vert _{2}+\left\Vert e^{-\nu
2^{\alpha q}t}\right\Vert _{L^{2}(\mathbb{R}^{+})}\left\Vert
F_{q}\right\Vert _{L^{\rho }([0,T])} \\
&\leq &\left( \frac{1-e^{-\nu 2^{\alpha q}T}}{\nu r}\right) ^{1/r}2^{-\frac{%
\alpha }{r}q}\left\Vert \Delta _{q}\theta _{\ast }\right\Vert _{2}+C_{\nu
}2^{-\frac{\alpha }{2}q}\left\Vert F_{q}\right\Vert _{L^{\rho }([0,T])}.
\end{eqnarray*}%
Multiplying the both sides of the last inequality by $2^{(s+\frac{%
\alpha }{r})q},$ using the commutator estimate (\ref{com}) and
then taking the $l^{2}(\mathbb{Z})$ norm, we obtain the desired
estimate (\ref{r2}).
\end{proof}

\begin{lemma}
\label{lem9}Let $\left( x_{n}\right) _{n\in \mathbb{N}}$ be a non negative
real sequence. Assume there exists two constants $A$ and $B\geq 0$ such that
\begin{equation*}
\left\{
\begin{array}{c}
x_{0}\leq 2A \\
\forall n,~x_{n+1}\leq A+Bx_{n}x_{n+1}.%
\end{array}%
\right.
\end{equation*}%
If $4AB\leq 1,$ then%
\begin{equation}
\forall n,~x_{n}\leq 2A.  \label{l}
\end{equation}
\end{lemma}

\begin{proof}
This lemma can be easily proved by induction.
\end{proof}

\begin{lemma}
\label{lem10}Let $\left( x_{n}\right) _{n\in \mathbb{N}}$ be a non negative
real sequence. Assume there exists two constants $A$ and $\delta \geq 0$
such that%
\begin{equation*}
\forall n,~x_{n+1}\leq A+\delta x_{n}.
\end{equation*}%
If $\delta <1,$ then the sequence $\left( x_{n}\right) _{n\in \mathbb{N}}$
is bounded.
\end{lemma}

\begin{proof}
Obvious. In fact, for all $n$ we have%
\begin{eqnarray*}
x_{n} &\leq &A\sum_{k=0}^{n}\delta ^{k}+\delta ^{n}x_{0} \\
&\leq &\frac{A}{1-\delta }+x_{0}.
\end{eqnarray*}
\end{proof}

\begin{lemma}
\label{lem11}let $\left( x_{n}\right) _{n\in \mathbb{N}}$ be a sequence in a
normed vectorial space $(E,\left\Vert .\right\Vert )$ Assume, there exists a
real sequence $\left( \gamma _{n,p}\right) _{n,p}$ and a real number $\delta
$ such that for all $(n,p)$ in $\mathbb{N}^{2},$%
\begin{equation}
\left\Vert x_{n+1+p}-x_{n+1}\right\Vert \leq \gamma _{n,p}+\delta \left\Vert
x_{n+p}-x_{n}\right\Vert .  \label{cau}
\end{equation}%
If%
\begin{equation*}
\delta \in \lbrack 0,1[\text{ and }\left( \sup_{p\in \mathbb{N}}\gamma
_{n,p}\right) _{n}\rightarrow 0\text{ as }n\rightarrow \infty
\end{equation*}%
then, the sequence $\left( x_{n}\right) _{n}$ is of Cauchy in $(E,\left\Vert
.\right\Vert ).$
\end{lemma}

\begin{proof}
Define $u_{n}=\sup_{p\in \mathbb{N}}\left\Vert x_{n+p}-x_{n}\right\Vert $
and $\gamma _{n}=\sup_{p\in \mathbb{N}}\gamma _{n,p}.$ According to (\ref%
{cau}), we have%
\begin{equation*}
\forall n,~u_{n+1}\leq \gamma _{n}+\delta u_{n}.
\end{equation*}%
This implies%
\begin{equation*}
\forall n,~u_{n}\leq \sum_{k=0}^{n}\delta ^{k}\gamma _{n-k}+\delta ^{n}u_{0}.
\end{equation*}%
Hence, for all $n>n_{0}$ in $\mathbb{N},$ we have%
\begin{equation*}
u_{n}\leq \frac{1}{1-\delta }\left( \sup_{0\leq k\leq n_{0}}\gamma
_{n-k}\right) +\left( \sup_{k\geq 0}\gamma _{k}\right) \frac{\delta
^{n_{0}+1}}{1-\delta }+\delta ^{n}u_{0}.
\end{equation*}%
Fixing $n_{0}$ and taking the $\overline{\lim_{n}},$ we obtain%
\begin{equation*}
\overline{\lim_{n}}~u_{n}\leq \left( \sup_{k\geq 0}\gamma _{k}\right) \frac{%
\delta ^{n_{0}+1}}{1-\delta }.
\end{equation*}%
Letting $n_{0}\rightarrow +\infty ,$ we get%
\begin{equation*}
\overline{\lim_{n}}~u_{n}=0.
\end{equation*}%
That is $\left( x_{n}\right) _{n}$ is a Cauchy sequence in $(E,\left\Vert
.\right\Vert ).$
\end{proof}

\section{Proof of Proposition \protect\ref{pr1}}

This section is devoted to the proof of Proposition \ref{pr1}. This proof is
motivated by the work \cite{CMZ}.

\begin{proof}[Proof of Proposition \protect\ref{pr1}]
Set $\theta _{0}\equiv 0.$ Proposition \ref{pr3} enables us to construct by
induction the sequence of functions $\theta _{n}\in C^{1}(\mathbb{R}%
^{+},H^{\infty })$ solutions to the systems%
\begin{equation*}
\left\{
\begin{array}{c}
\partial _{t}\theta _{n+1}+\Lambda ^{\alpha }\theta _{n+1}+u_{n}\vec{\nabla}%
\theta _{n+1}=0 \\
\left. \theta _{n+1}\right\vert _{t=0}=S_{n+1}\theta _{\ast } \\
u_{n}=\Lambda ^{\alpha -1}\mathcal{R}^{\perp }(\theta _{n})%
\end{array}%
\right.
\end{equation*}%
and satisfying, for all $r\geq 2,~s\geq 0$ and $T>0,$ the following
estimates:%
\begin{equation}
\sup_{t\geq 0}\left\Vert \theta _{n}(t)\right\Vert _{2}\leq \left\Vert
\theta _{\ast }\right\Vert _{2}  \label{est1}
\end{equation}%
\begin{equation}
\left\Vert \theta _{n+1}\right\Vert _{\tilde{L}_{T}^{r}\dot{H}^{s+\frac{%
\alpha }{r}}}\leq \mathcal{K}_{r,s}(\theta _{\ast },T)+\mathcal{C}%
_{r,s}\left( \left\Vert \theta _{n}\right\Vert _{L_{T}^{2}\dot{H}^{1+\frac{%
\alpha }{2}}}\left\Vert \theta _{n+1}\right\Vert _{\tilde{L}_{T}^{r}\dot{H}%
^{s+\frac{\alpha }{r}}}+\left\Vert \theta _{n+1}\right\Vert _{L_{T}^{2}\dot{H%
}^{1+\frac{\alpha }{2}}}\left\Vert \theta _{n}\right\Vert _{\tilde{L}_{T}^{r}%
\dot{H}^{s+\frac{\alpha }{r}}}\right)  \label{est2}
\end{equation}%
where $\mathcal{K}_{r,s}(\theta _{\ast },T)$ is defined by (\ref{r3})-(\ref%
{inf}) and $\mathcal{C}_{r,s}$ is a constant depending only on $r$ and $s.$

Let $\varepsilon _{\sigma }>0$ to be determined later and assume that for
some $T>0$ we have
\begin{equation*}
\mathcal{K}(\theta _{\ast },T)=\mathcal{K}_{2,1}(\theta _{\ast },T)\leq
\varepsilon _{\sigma },
\end{equation*}%
($T$ exists since $\mathcal{K}(\theta _{\ast },T)\rightarrow 0$ as $%
T\rightarrow 0^{+}$). Then, applying the estimate (\ref{est2}) with $s=1$
and $r=2$ yields%
\begin{equation*}
\left\Vert \theta _{n+1}\right\Vert _{L_{T}^{2}\dot{H}^{1+\frac{\alpha }{2}%
}}\leq \varepsilon _{\sigma }+2\mathcal{C}_{2,1}\left\Vert \theta
_{n}\right\Vert _{L_{T}^{2}\dot{H}^{1+\frac{\alpha }{2}}}\left\Vert \theta
_{n+1}\right\Vert _{L_{T}^{2}\dot{H}^{1+\frac{\alpha }{2}}}.
\end{equation*}%
Hence, according to Lemma \ref{lem9}, we have%
\begin{equation}
\forall n,~\left\Vert \theta _{n}\right\Vert _{L_{T}^{2}\dot{H}^{1+\frac{%
\alpha }{2}}}\leq 2\varepsilon _{\sigma }  \label{est3}
\end{equation}%
provided that%
\begin{equation}
8\varepsilon _{\sigma }\mathcal{C}_{2,1}\leq 1.  \label{Cond1}
\end{equation}%
Using now the estimate (\ref{est2}) with $s=\sigma $ and $r\in \{2,+\infty
\},$ we get%
\begin{equation*}
\left\Vert \theta _{n+1}\right\Vert _{\tilde{L}_{T}^{r}\dot{H}^{\sigma +%
\frac{\alpha }{r}}}\leq \mathcal{K}_{r,s}(\theta _{\ast },T)+2\varepsilon
_{\sigma }\mathcal{C}_{r,\sigma }\left( \left\Vert \theta _{n}\right\Vert _{%
\tilde{L}_{T}^{r}\dot{H}^{\sigma +\frac{\alpha }{r}}}+\left\Vert \theta
_{n+1}\right\Vert _{\tilde{L}_{T}^{r}\dot{H}^{\sigma +\frac{\alpha }{r}%
}}\right) .
\end{equation*}%
Therefore, Lemma \ref{lem10} ensures that the sequence $\left(
\theta _{n}\right) _{n}$ is bounded in the spaces
$\tilde{L}_{T}^{2}\dot{H}^{\sigma +\frac{\alpha }{2}}$ and
$\tilde{L}_{T}^{\infty }\dot{H}^{\sigma }$ provided
that%
\begin{equation}
4\varepsilon _{\sigma }\max \{\mathcal{C}_{2,\sigma },\mathcal{C}_{\infty
,\sigma }\}\leq 1.  \label{Cond2}
\end{equation}%
We claim now that if $\varepsilon _{\sigma }$ is small enough then the
sequence $\left( \theta _{n}\right) _{n}$ is of Cauchy in the space $%
C([0,T],L^{2}(\mathbb{R}^{2})).$ Let $n\in \mathbb{N}$ and $p\in \mathbb{N}%
^{\ast }.$ Define $\omega _{n+1}=\theta _{n+1+p}-\theta _{n+1}$, $\omega
_{n}=\theta _{n+p}-\theta _{n}$ and $v_{n}=\Lambda ^{\alpha -1}\mathcal{R}%
^{\perp }(\omega _{n}).$ Clearly, we have the equation%
\begin{equation*}
\partial _{t}\omega _{n+1}+\Lambda ^{\alpha }\omega _{n+1}+u_{n+p}\vec{\nabla%
}\omega _{n+1}+v_{n}\vec{\nabla}\theta _{n+1}=0.
\end{equation*}%
Hence, Lemma \ref{lem5} implies that for all $t$ in $[0,T],$
\begin{equation}
\frac{d}{2dt}\left\Vert \omega _{n+1}(t)\right\Vert _{2}^{2}+\left\Vert
\Lambda ^{\alpha /2}\omega _{n+1}(t)\right\Vert _{2}^{2}\leq \kappa
_{n}(t)\equiv \left\vert \langle v_{n}\vec{\nabla}\theta _{n+1},\omega
_{n+1}\rangle \right\vert .  \label{est4}
\end{equation}%
Now a simple calculation gives the following estimates%
\begin{eqnarray}
\kappa _{n}(t) &\leq &\left\Vert v_{n}\vec{\nabla}\theta _{n+1}\right\Vert _{%
\dot{H}^{-\alpha /2}}\left\Vert \omega _{n+1}\right\Vert _{\dot{H}^{\alpha
/2}}  \label{l1} \\
&\lesssim &\left\Vert v_{n}(t)\right\Vert _{\dot{H}^{1-\alpha }}\left\Vert
\vec{\nabla}\theta _{n+1}(t)\right\Vert _{\dot{H}^{\alpha /2}}\left\Vert
\Lambda ^{\alpha /2}\omega _{n+1}(t)\right\Vert _{2}  \label{l2} \\
&\lesssim &\left\Vert \omega _{n}(t)\right\Vert _{2}\left\Vert \theta
_{n+1}(t)\right\Vert _{\dot{H}^{1+\alpha /2}}\left\Vert \Lambda ^{\alpha
/2}\omega _{n+1}(t)\right\Vert _{2}  \notag \\
&\leq &C\left\Vert \theta _{n+1}(t)\right\Vert _{\dot{H}^{1+\alpha
/2}}^{2}\left\Vert \omega _{n}(t)\right\Vert _{2}^{2}+\left\Vert \Lambda
^{\alpha /2}\omega _{n+1}(t)\right\Vert _{2}^{2}  \label{l4}
\end{eqnarray}%
where to obtain (\ref{l2}) from (\ref{l1}) we have used Lemma (\ref{lem3})
with $\sigma _{1}=1-\alpha $ and $\sigma _{2}=\alpha /2.$

Substituting (\ref{l4}) in (\ref{est4}), integrating with respect to time
and taking the supremum over all $t\in \lbrack 0,T]$ yield
\begin{equation*}
\sup_{0\leq t\leq T}\left\Vert \omega _{n+1}(t)\right\Vert _{2}^{2}\leq
\left\Vert \omega _{n+1}(0)\right\Vert _{2}^{2}+C\sup_{0\leq t\leq
T}\left\Vert \omega _{n+1}(t)\right\Vert _{2}^{2}\int_{0}^{T}\left\Vert
\theta _{n+1}(t)\right\Vert _{\dot{H}^{1+\alpha /2}}^{2}dt.
\end{equation*}%
Recalling the estimate (\ref{est3}) and using the the fact%
\begin{equation*}
\forall a,b\in \mathbb{R}^{+},~\sqrt{a^{2}+b^{2}}\leq a+b,
\end{equation*}%
we infer from the above inequality that%
\begin{equation}
\left\Vert \theta _{n+1+p}-\theta _{n+1}\right\Vert _{\ast }\leq \delta
_{n,p}+C\varepsilon _{\sigma }\left\Vert \theta _{n+p}-\theta
_{n}\right\Vert _{\ast }  \label{est5}
\end{equation}%
where $\left\Vert .\right\Vert _{\ast }$ denotes the norm of the space $%
L^{\infty }([0,T],L^{2}(\mathbb{R}^{2}))$ and
\begin{equation*}
\delta _{n,p}\equiv \left\Vert S_{n+p+1}\theta _{\ast }-S_{n+1}\theta _{\ast
}\right\Vert _{2}.
\end{equation*}%
Now since $\theta _{\ast }\in L^{2}(\mathbb{R}^{2})$ then $\left(
\sup_{p}\delta _{n,p}\right) _{n}\rightarrow 0$ as $n\rightarrow +\infty .$
Hence, according to Lemma \ref{lem11}, the sequence $\left( \theta
_{n}\right) _{n}$ is of Cauchy in the space $C([0,T],L^{2}(\mathbb{R}^{2}))$
provided%
\begin{equation}
C\varepsilon _{\sigma }<1.  \label{Cond3}
\end{equation}%
In conclusion, there exists a constant $\varepsilon _{\sigma }>0$ depending
only on $\sigma $ such that if $\mathcal{K}(\theta _{\ast },T)\leq
\varepsilon _{\sigma }$ for some $T>0$ then the sequence $(\theta _{n})_{n}$
is bounded in the space~$\mathbf{X}_{T}^{\sigma }\equiv \tilde{L}%
_{T}^{\infty }H^{\sigma }\cap L_{T}^{2}H^{\sigma +\frac{\alpha }{2}}$ and
converges in the space $C([0,T],L^{2}(\mathbb{R}^{2}))$ to a function $%
\theta .$ Hence, according to Lemma \ref{lem6} and Lemma \ref{lem7}, the
function $\theta $ belongs to the space $\mathbf{X}_{T}^{\sigma }\cap
C([0,T],H^{\sigma }(\mathbb{R}^{2})).$ On the other hand, using the
embedding $\tilde{L}_{T}^{\infty }H^{\sigma }\subset L_{T}^{\infty
}H^{\sigma }$ and the interpolation inequality%
\begin{equation*}
\forall 0\leq s\leq \sigma ,~\left\Vert f\right\Vert _{\dot{H}^{s}}\leq
\left\Vert f\right\Vert _{2}^{1-\frac{s}{\sigma }}\left\Vert f\right\Vert _{%
\dot{H}^{\sigma }}^{\frac{s}{\sigma }},
\end{equation*}%
we deduce that $\left( \theta _{n}\right) _{n}$ converges to $\theta $ in $%
C([0,T],\dot{H}^{\sigma ^{\prime }}(\mathbb{R}^{2}))$ for all $\sigma
^{\prime }<\sigma $ which implies in particular that
\begin{eqnarray*}
\Lambda ^{\alpha }\theta _{n} &\rightarrow &\Lambda ^{\alpha }\theta \text{
in }C([0,T],L^{2}(\mathbb{R}^{2})) \\
\theta _{n} &\rightarrow &\theta \text{ in }C([0,T],\dot{H}^{\sigma _{\ast
}}(\mathbb{R}^{2})) \\
u_{n} &\rightarrow &u\equiv \Lambda ^{\alpha -1}\mathcal{R}^{\perp }(\theta )%
\text{ in }C([0,T],\dot{H}^{\sigma _{\ast }+1-\alpha }(\mathbb{R}^{2}))
\end{eqnarray*}%
where $\sigma _{\ast }=\frac{1+\alpha }{2}.$ Thus, using the fact that $u_{n}%
\vec{\nabla}\theta _{n+1}=\nabla .(\theta _{n+1}u_{n})$ and Lemma \ref{lem3}
with $\sigma _{1}=\sigma _{\ast }+1-\alpha $ and $\sigma _{2}=\sigma _{\ast
},$ we deduce that $u_{n}\vec{\nabla}\theta _{n+1}$ converges to $u\vec{%
\nabla}\theta $ in $C([0,T],L^{2}(\mathbb{R}^{2})).$ We then
conclude that the function $\theta $ belongs to the space
$\mathbf{Z}_{T}^{\sigma }$ and
satisfies%
\begin{equation*}
\left\{
\begin{array}{c}
\partial _{t}\theta +\Lambda ^{\alpha }\theta +u\vec{\nabla}\theta =0 \\
u=\Lambda ^{\alpha -1}\mathcal{R}^{\perp }(\theta ) \\
\left. \theta \right\vert _{t=0}=\theta _{\ast }.%
\end{array}%
\right.
\end{equation*}%
Finally, it remains to prove the uniqueness. Assume that $\theta _{a}$ and $%
\theta _{b}\in \mathbf{Z}_{T}^{\sigma }$ are two solutions to the equation
(MQG) with the same data $\theta _{\ast }.$ Set $u_{a}=\Lambda ^{\alpha -1}%
\mathcal{R}^{\perp }(\theta _{a}),~u_{b}=\Lambda ^{\alpha -1}\mathcal{R}%
^{\perp }(\theta _{b})$, $\omega =\theta _{a}-\theta _{b}$ and $%
u=u_{a}-u_{b}.$ We have the equation%
\begin{equation*}
\partial _{t}\omega +\Lambda ^{\alpha }\omega +u_{a}\vec{\nabla}\omega +u%
\vec{\nabla}\theta _{b}=0
\end{equation*}%
which implies, by Lemma \ref{lem5},%
\begin{equation}
\frac{1}{2}\frac{d}{dt}\left\Vert \omega (t)\right\Vert _{2}^{2}+\left\Vert
\Lambda ^{\alpha /2}\omega (t)\right\Vert _{2}^{2}\leq \kappa (t)\equiv
\left\vert \int_{\mathbb{R}^{2}}u\vec{\nabla}\theta _{b}\omega dx\right\vert
.  \label{eq}
\end{equation}%
Following the same argument leading to (\ref{l4}), we obtain%
\begin{equation*}
\kappa (t)\leq C\left\Vert \theta _{b}(t)\right\Vert _{\dot{H}^{1+\alpha
/2}}^{2}\left\Vert \omega (t)\right\Vert _{2}^{2}+\left\Vert \Lambda
^{\alpha /2}\omega (t)\right\Vert _{2}^{2}.
\end{equation*}%
Inserting this estimate in the inequality (\ref{eq}), we get%
\begin{equation*}
\frac{d}{dt}\left\Vert \omega (t)\right\Vert _{2}^{2}\lesssim \left\Vert
\theta _{b}(t)\right\Vert _{\dot{H}^{1+\alpha /2}}^{2}\left\Vert \omega
(t)\right\Vert _{2}^{2}.
\end{equation*}%
Thereby, the Gronwall inequality implies $\omega =0$ that is $\theta
_{a}=\theta _{b}.$ This completes the proof of the proposition.
\end{proof}

\section{Proof of Proposition \protect\ref{pr2}}

The proof of Proposition \ref{pr2} is based on the following blowup
criterion.

\begin{lemma}
\label{lem12}Let $\theta _{\ast }\in H^{\sigma }(\mathbb{R}^{2})$ with $%
\sigma \geq 1$ and let
\begin{equation*}
\theta \in \bigcap_{T<T^{\ast }}\mathbf{Z}_{T}^{\sigma }
\end{equation*}%
be the maximal solution to the equation (MQG) with initial data $\theta
_{\ast }.$ Assume $T^{\ast }<\infty .$ Then there is a constant $c_{\sigma
}>0$ depending only on $\sigma $ such that
\begin{equation}
\forall 0\leq t<T^{\ast },~\left\Vert \theta (t)\right\Vert _{\dot{H}^{1+%
\frac{\alpha }{2}}}\geq \frac{c_{\sigma }}{\sqrt{T^{\ast }-t}}.  \label{eg}
\end{equation}%
In particular,%
\begin{equation*}
\int_{0}^{T^{\ast }}\left\Vert \theta (t)\right\Vert _{\dot{H}^{1+\frac{%
\alpha }{2}}}^{2}dt=+\infty .
\end{equation*}
\end{lemma}

\begin{proof}
Let $t$ in $[0,T^{\ast }[.$ Firstly, according to the last assertion of
Proposition \ref{pr1}, we must have
\begin{equation}
\mathcal{K}(\theta (t),T^{\ast }-t)\geq \varepsilon _{\sigma }.
\end{equation}%
Secondly, for every $f\in H^{1+\frac{\alpha }{2}}$ and $T>0,$ we have the
following estimates%
\begin{eqnarray}
\mathcal{K}(f,T) &\leq &\sqrt{T}\sup_{q\in \mathbb{Z}}\left( \frac{%
1-e^{-2\nu 2^{\alpha q}T}}{2\nu 2^{\alpha q}T}\right) ^{1/2}\left\Vert
\left( 2^{(1+\frac{\alpha }{2})q}\left\Vert \Delta _{q}f\right\Vert
_{2}\right) _{q}\right\Vert _{l^{2}(\mathbb{Z})}  \notag \\
&\leq &C_{\ast }\sqrt{T}\left\Vert f\right\Vert _{\dot{H}^{1+\frac{\alpha }{2%
}}}
\end{eqnarray}%
where $C_{\ast }=\sup_{x>0}\sqrt{\frac{1-e^{-x}}{x}}.$ Hence, combining the
estimates (\ref{eg1})-(\ref{eg2}) yields (\ref{eg}) with $c_{\sigma }=\frac{%
\varepsilon _{\sigma }}{C_{\ast }}.$
\end{proof}

\bigskip

Now we are in position to prove our proposition:

\begin{proof}[Proof of Proposition \protect\ref{pr2}]
In view of Proposition \ref{pr1}, the modified quasi-geostrophic equation
corresponding to the initial data $\tilde{\theta}_{\ast }\equiv \theta
(t_{0})$ has a unique maximal solution
\begin{equation*}
\tilde{\theta}\in \cap _{T<T^{\ast }}\mathbf{Z}_{T}^{\sigma }.
\end{equation*}%
On the other hand, one can easily verify that the function $\theta
_{t_{0}}\equiv \theta (.+t_{0})$ is also a solution to the same (MQG)
equation. Hence uniqueness in the space $\mathbf{Z}_{T}^{1}$ implies that $%
\theta _{t_{0}}=\tilde{\theta}$ on the interval $[0,\tau _{\ast }[$ where $%
\tau _{\ast }=\inf \{T^{\ast };T_{\ast }-t_{0}\}.$ Therefore, since $\theta
_{t_{0}}$ belongs to the space $L^{2}([0,T_{\ast }-t_{0}[,H^{1+\frac{\alpha
}{2}})$ then the preceding lemma ensures that $T^{\ast }>T_{\ast }-t_{0}.$
Thus we deduce that $\theta _{t_{0}}\in \mathbf{Z}_{T_{\ast }-t_{0}}^{\sigma
}$ which implies the required result.
\end{proof}

\section{Proof of the main Theorem}

Firstly, Proposition \ref{pr1} ensures the existence of a unique maximal
solution $\theta \in \cap _{T<T^{\ast }}\mathbf{Z}_{T}^{\sigma }$ to the
equation (MQG). Now let $a<T$ be a two fixed real-number in the interval $%
]0,T^{\ast }[$. successive application of Proposition 2, allows us to
construct an increasing sequence $(t_{n})_{n\in \mathbb{N}}\in ]0,a[$ such
that for all $n,$ $\theta (t_{n})\in H^{\sigma _{n}}$ and
\begin{equation*}
\theta \in C([t_{n},T],H^{\sigma _{n}})\cap L^{2}([t_{n},T];H^{\sigma
_{n+1}})
\end{equation*}%
where $\sigma _{k}\equiv \sigma +k\frac{\alpha }{2}.$ Consequently, the
solution $\theta $ belongs to the space $C([a,T],H^{\infty }).$ Now, since $%
a $ and $T$ are arbitrary chosen in $]0,T^{\ast }[$ then $\theta $ is in the
space $C\left( ]0,T^{\ast }[,H^{\infty }\right) .$ On the other hand, the
equation
\begin{equation*}
\partial _{t}\theta =-\Lambda ^{\alpha }\theta -u\vec{\nabla}\theta ,
\end{equation*}%
combined with the continuity $\Lambda ^{\alpha }$ and $\vec{\nabla}$ on the
space $H^{\infty }$ and Lemma \ref{lem1}, enable us, via a standard
Boot-strap argument, to convert the space regularity of $\theta $ to time
regularity. We then deduce that
\begin{equation*}
\theta \in C^{\infty }\left( ]0,T^{\ast }[,H^{\infty }\right) .
\end{equation*}%
Now, we will establish that the solution $\theta $ is global in time, that
is $T^{\ast }=\infty .$ We will argue by opposition: we suppose that $%
T^{\ast }<\infty .$ First, from Lemma \ref{lem5}, we easily see that $\theta
$ is a Leray-Hopf solution to the (MQG) equation, that is%
\begin{equation*}
\theta \in L_{T^{\ast }}^{\infty }L^{2}\cap L_{T^{\ast }}^{2}\dot{H}^{\frac{%
\alpha }{2}}.
\end{equation*}%
Hence, from the papers \cite{CIW} and \cite{CW2} we deduce that for all $%
t_{0}$ in $]0,T^{\ast }[$
\begin{equation*}
\sup_{t_{0}\leq t<T^{\ast }}\left\Vert \theta (t)\right\Vert _{C_{B}^{2}(%
\mathbb{R}^{2})}<\infty .
\end{equation*}%
Therefore, the last assertion of Lemma \ref{lem2} implies%
\begin{equation*}
{\footnotesize M}_{t_{0}}\equiv \sup_{t_{0}\leq t<T^{\ast }}\left\Vert
u(t)\right\Vert _{C_{B}^{2}(\mathbb{R}^{2})}<\infty .
\end{equation*}%
Fix $t_{0}$ in $]0,T^{\ast }[$ and let $\beta \in \mathbb{N}^{2}$ any
multi-index with $\left\vert \beta \right\vert \leq 2.$ The function $%
D^{\beta }\theta $ satisfies%
\begin{equation*}
\partial _{t}D^{\beta }\theta +\Lambda ^{\alpha }D^{\beta }\theta +u\vec{%
\nabla}D^{\beta }\theta =\left[ u,D^{\beta }\right] \vec{\nabla}\theta .
\end{equation*}%
Thus, Lemma \ref{lem5} implies that for all $t$ in $[t_{0},T^{\ast }[$%
\begin{eqnarray*}
\frac{d}{dt}\left\Vert D^{\beta }\theta (t)\right\Vert _{2}^{2} &\lesssim
&\left\Vert \left[ u,D^{\beta }\right] \vec{\nabla}\theta (t)\right\Vert
_{2}\left\Vert D^{\beta }\theta (t)\right\Vert _{2} \\
&\lesssim &{\footnotesize M}_{t_{0}}\left\Vert \theta (t)\right\Vert
_{H^{2}}^{2}.
\end{eqnarray*}%
Summing on $\beta ,$ we get for all $t$ in $[t_{0},T^{\ast }[$%
\begin{equation*}
\frac{d}{dt}\left\Vert \theta (t)\right\Vert _{H^{2}}^{2}\lesssim
{\footnotesize M}_{t_{0}}\left\Vert \theta (t)\right\Vert _{H^{2}}^{2}.
\end{equation*}%
Thanks to The Gronwall inequality, this inequation implies that
\begin{equation*}
\sup_{t_{0}\leq t<T^{\ast }}\left\Vert \theta (t)\right\Vert _{H^{2}}<\infty
\end{equation*}%
which contradicts Lemma \ref{lem12}, since $H^{2}\hookrightarrow \dot{H}^{1+%
\frac{\alpha }{2}}$. We then conclude that $T^{\ast }=\infty .$

\bigskip

\providecommand{\bysame}{\leavevmode\hbox
to3em{\hrulefill}\thinspace}



\providecommand{\bysame}{\leavevmode\hbox
to3em{\hrulefill}\thinspace} \providecommand{\MR}{\relax\ifhmode\unskip%
\space\fi MR }
\providecommand{\MRhref}[2]{  \href{http://www.ams.org/mathscinet-getitem?mr=#1}{#2}
} \providecommand{\href}[2]{#2}

\end{document}